\documentclass[a4paper]{article}
\usepackage{latexsym}

\reversemarginpar
\newcommand{\kitem}{\begin{itemize}\vspace{-2ex}}
\newcommand{\kenditem}{\vspace{-1ex}\end{itemize}}

\newcommand{\Z}{{Z\!\!\!Z}}
\newcommand{\Q}{I\!\!\!\!Q}
\newcommand{\R}{{I\!\!R}}
\newcommand{\C}{{\,I\!\!\!\!C}}
\renewcommand{\L}{I\!\!L}
\renewcommand{\P}{{I\!\!P}}
\newcommand{\CO}{{\cal O}}

\renewcommand{\subseteq}{\subset}

\newcommand{\ko}{\overline}
\newcommand{\ku}{\underline}
\newcommand{\ks}{\scriptstyle}
\newcommand{\kss}{\scriptscriptstyle}
\newcommand{\kd}{\displaystyle}
\newcommand{\kf}{\footnotesize}

\def\mapup#1{\llap{$\vcenter{\hbox{$\scriptstyle #1$}}$}\Big\uparrow}
\def\cdmatrix#1{\def\normalbaselines{\baselineskip20pt
 \lineskip3pt \lineskiplimit3pt }
\matrix{#1}}

\newcommand{\qed}{{
\unskip\nobreak\hfil\penalty50\hskip0.1em\hbox{}\nobreak\hfil$\Box$
\parfillskip=0pt \finalhyphendemerits=0 \par\medbreak}}

\newcommand{\kkk}[1]{\marginpar{\raggedright\sc\small #1}}
\def\kkk#1{}

\def\thesection{\arabic{section}.}
\def\thesubsection{\thesection\arabic{subsection}}

\newcommand{\neu}[1]{\protect\refstepcounter{subsection}\protect
   \label{#1}\vspace{1ex}
   {\bf (\thesubsection)} $\quad$\kkk{#1}\ignorespaces}
\def\zitat#1#2{(\ref{#1-#2})} 

\newcommand{\kb}{{\kss \bullet}}
\renewcommand{\span}{\mbox{\rm span}}
\newcommand{\im}{\mbox{\rm im}\,}

\newcommand{\rank}{\mbox{\rm rank}}
\newcommand{\innt}{\mbox{\rm int}\,}
\newcommand{\cg}{{\rm cg}}
\newcommand{\scg}{\cg}
\newcommand{\gHom}{\mbox{\rm Hom}}
\newcommand{\sh}{{\rm sh}}
\newcommand{\kn}{{n}}
\newcommand{\Harr}{\mathop{\rm Harr}\nolimits}
\newcommand{\HA}{H\!A}
\newcommand{\m}{\ku{m}}
\newcommand{\kht}{{\rm ht}}
\newcommand{\sht}{\kht}
\begin{document}
\title{\bf Cotangent cohomology of rational  surface singularities}
\author{\Large Klaus Altmann \qquad\qquad Jan Stevens}
\date{}
\maketitle

\begin{abstract}
In this paper we show that the number of generators of
the cotangent cohomology groups $T_Y^\kn$, $\kn\geq2$, is the
same for  all rational  surface singularities $Y$.
For  a large class of rational  surface singularities, including
quotient singularities, this number is also the dimension. 
For them we
obtain an explicit  formula for the Poincar\'{e} series $P_Y(t)
=\sum \dim\,T^\kn_Y\cdot t^n$. In the special case 
of the cone over the rational normal curve we give the multigraded
Poincar\'e series.  

\end{abstract}

%%%%%%%%%%%%%%%%%%
%
% Introduction
%
%%%%%%%%%%%%%%%%%%
\section{Introduction}\label{intro}

The cotangent cohomology groups $T^\kn$ with small $\kn$ play an important role
in the deformation theory of singularities: 
$T^1$ classifies infinitesimal deformations   and the obstructions
land in $T^2$. Originally constructed ad hoc, the correct way to obtain 
these groups is as the cohomology of the cotangent complex. This yields 
also higher $T^\kn$, which no longer have a direct meaning in terms
of deformations.
\par

In this paper we study these higher 
cohomology groups $T_Y^\kn$ for   rational  surface singularities $Y$.
For  a large class of rational  surface singularities, including
quotient singularities, we obtain their dimension. 
For an explicit formula for the Poincar\'{e} series $P_Y(t)$, 
see \zitat{hps}{app}.
\par

Our methods are a combination of the following three items:
\kitem
\item[(1)]
We use the hyperplane section machinery of \cite{BC}
to move freely between surface singularities, partition curves, 
and  fat points. 
It  suffices to compute the cohomology groups 
$T_Y^\kn$ for special singularities, to obtain the number of
generators for all rational surface singularities.
\item[(2)]
In many cases, cotangent cohomology may be obtained via Harrison cohomology,
which is much easier to handle.
Using a Noether normalisation the Harrison complex gets linear 
over a bigger ring than just $\C$ (which is our ground field).
\item[(3)]
Taking for $Y$ a cone over a rational normal
curve, we may use the explicit description 
of $T_Y^\kn$  obtained in \cite{AQ} by toric
methods. 
\kenditem
The descriptions in (2) and (3) complement each other and
show that $T^\kn_Y$ of the cone over the
rational normal curve  is concentrated in degree $-\kn$.
This allows us to compute the dimension as Euler characteristic.

The paper is organised as follows.
After recalling the definitions of cotangent and Harrison cohomology
we review its computation for the case of
the fat point of minimal multiplicity
and give the explicit formula for its
Poincar\'e series (we are indebted to Duco van Straten and
Ragnar Buchweitz for help on this point).\\
Section 3 describes the applications of Noether normalisation to
the computation of Harrison cohomology. The main result is the
degree bound for the cotangent cohomology of Cohen-Macaulay singularities
of minimal multiplicity from below, cf.\ Corollary \zitat{noeth}{zero}(2).\\
In the next section toric methods are used to deal with the cone over the
rational normal curve. In this special case we can bound the degree of the
cohomology groups from above, too. As a consequence, we obtain complete
information about the Poincar\'e series.\\
Finally, using these results  as input for the hyperplane
machinery we find in the last section the Poincar\'e series for the
partition curves and obtain that their $T^\kn$ is
annihilated by the maximal ideal. This then implies that the
number of generators of $T^\kn_Y$ is the same for all rational
surface singularities.
\par

{\bf Notation:}
We would like to give the following 
guide line concerning the notation for Poincar\'e series.
The symbol $Q$ denotes those series involving Harrison cohomology of
the actual space or ring with values in $\C$, 
while  $P$ always points to
the usual cotangent cohomology of the space  itself.
Moreover, if these letters come with a tilde, then a finer grading than
the usual $\Z$-grading is involved.
\par

%%%%%%%%%%
%
%  Cotangent cohomology and Harrison cohomology
%
%%%%%%%%%
\section{Cotangent cohomology and Harrison cohomology}\label{cotan}

%%%%%%%%%%%%%%%
% (cotan.gen)
%%%%%%%%%%%%%%%

\neu{cotan-gen}
Let $A$ be a commutative algebra of essentially finite type
over a base-ring $S$. For any $A$-module $M$, one
gets the {\it Andr\'e--Quillen\/} or {\it cotangent cohomology groups\/} as
\[
T^\kn(A/S,M):=H^\kn\big({\rm Hom}_A(\L^{A/S}_*,M)\big)
\]
with $\L^{A/S}_*$ being the so-called cotangent complex.
We are going to recall the major properties of this cohomology theory. For
the details, including the definition of $\L^{A/S}_*$,
see \cite{Loday}.
\par

If $A$ is a smooth $S$-algebra, then $T^\kn(A/S,M)=0$ for $\kn\geq1$
and all $A$-modules $M$. 
For general $A$, a short exact sequence of $A$-modules gives a long
exact sequence in cotangent cohomology. Moreover, the Zariski-Jacobi sequence
takes care of ring homomorphisms $S \to A \to B$; for a $B$-module $M$
it looks like
\[
\cdots \longrightarrow T^\kn(B/A,M)
\longrightarrow T^\kn(B/S,M)
\longrightarrow T^\kn(A/S,M)
\longrightarrow T^{\kn+1}(B/A,M)
\longrightarrow \cdots
\]
The cotangent cohomology behaves well under base change.
Given a co-cartesian diagram
\[
\cdmatrix{
A&\longrightarrow &A' \cr
\mapup{\phi}&&\mapup{}\cr
S&\longrightarrow &S' \cr}
\]
with $\phi$ flat, and  an $A'$-module $M'$, there is a natural isomorphism
\[
T^\kn(A'/S',M')\cong T^\kn(A/S,M')\,.
\]
If, moreover, $S'$ is a flat $S$-module, then for any $A$-module $M$
\[
T^\kn(A'/S',M\otimes_S S')\cong T^\kn(A/S,M)\otimes_S S'\,.
\]
\par

%%%%%%%%%%%%%%%
% (cotan.harr)
%%%%%%%%%%%%%%%

\neu{cotan-harr}
To describe Harrison cohomology, we first recall Hochschild cohomology.
While this concept works also for non-commutative unital algebras, we
assume here the same setting as before.
For an $A$-module $M$, we consider the complex
\[
C^\kn(A/S,M):=\mbox{\rm Hom}_S(A^{\otimes \kn}, M)
\]
with differential 
\[
\displaylines{\qquad
(\delta f) (a_0,\dots,a_\kn):=
{}\hfill\cr\hfill
 a_0f(a_1,\dots,a_\kn)+
\sum_{i=1}^{\kn}(-1)^i f(a_0,\dots,a_{i-1}a_{i},\dots,a_\kn)
+(-1)^{\kn+1}a_\kn f(a_0,\dots,a_{\kn-1})\,.
\qquad}
\]
{\em Hochschild cohomology} $HH^\kn(A/S,M)$ is the cohomology of this complex. 
It can also
be computed from the so-called {\em reduced} subcomplex $\ko{C}^\kb(A/S,M)$
consisting only of those maps $f\colon A^{\otimes \kn}\to M$ that vanish 
whenever at least one of the arguments equals 1.
\par

{\bf Definition:}
A permutation $\sigma\in S_\kn$ is called a $(p,\kn-p)$-shuffle 
if $\,\sigma(1)<\dots<\sigma(p)$ and
$\,\sigma(p+1)<\dots<\sigma(\kn)$. Moreover, in the group algebra
$\Z[S_\kn]$ we define the elements
\vspace{-2ex}
\[
\sh_{p,\kn-p}:=\sum_{\mbox{\kf $(p,\kn-p)$-shuffles}} 
\hspace{-1.5em}\mbox{sgn}(\sigma)\, \sigma 
\hspace{2em}\mbox{and}\hspace{2em}
\sh:=\sum_{p=1}^{\kn-1}\mbox{\rm sh}_{p,\kn-p}\,.
\vspace{-1ex}
\]
\par

The latter element $\,\sh\in\Z[S_\kn]$ gives rise to the 
so-called shuffle invariant subcomplexes
\[
C_\sh^\kn(A/S,M):=
\big\{f\in \gHom_S(A^{\otimes \kn}, M)\bigm|
f(\sh(\ku{a}))=0 \; \mbox{ for every }\; \ku{a}\in A^{\otimes \kn}
\big\}
\subseteq C^\kn(A/S,M)
\]
and $\ko{C}_\sh^\kn(A/S,M)\subseteq \ko{C}^\kn(A/S,M)$
defined in the same manner. Both complexes yield the same cohomology, which is 
called  {\em Harrison cohomology\/}: 
\[
\Harr^\kn(A/S,M):= H^\kn\big( C_\sh^\kb(A/S,M)\big)
               = H^\kn\big( \ko{C}_\sh^\kb(A/S,M)\big)\,.
\vspace{-2ex}
\]
\par

%%%%%%%%%%%%%%%
% (cotan.compare)
%%%%%%%%%%%%%%%

\neu{cotan-compare}
The following well known result compares the cohomology theories defined so
far. Good references are \cite{Loday} or \cite{Pal}.
\par

{\bf Theorem:}
{\em
If\/ $\Q\subseteq S$, then
Harrison cohomology is a direct summand of Hochschild cohomology.
Moreover, if $A$ is a {\it flat\/} $S$-module, then
\[
\,T^{\kn}(A/S,M) \cong \Harr^{\kn+1}(A/S,M)\,.
\vspace{-3ex}
\]
}
\par

%%%%%%%%%%%%%%%
% (cotan.example)
%%%%%%%%%%%%%%%

\neu{cotan-example}
As an example, we consider the fat point $Z_m$ ($m\geq 2$)
with minimal multiplicity $d=m+1$.
Let $V$ be an $m$-dimensional $\C$-vector space and let 
$A=\CO_{Z_m}$ be the ring 
$\C\oplus V$ with trivial multiplication $V^2=0$. 
First we compute
the Hochschild cohomology $HH^\kb (A/\C,A)$.
The reduced complex is
\[
\ko C^\kn(A/\C,A) = \gHom_\C(V^{\otimes \kn},A)\,.
\]
Because $ab=0\in A$ for all $a,b\in V$, the differential reduces to 
\[
(\delta f) (a_0,\dots,a_\kn)= a_0f(a_1,\dots,a_\kn)+
(-1)^{\kn+1}a_\kn f(a_0,\dots,a_{\kn-1})\,.
\]
We conclude that $\delta f =0 $ if and only if $\;\im f \subset V$;
hence
\[
HH^\kn(A/\C,A)= \gHom\,(V^{\otimes \kn},V)\Big/
\delta \,\gHom\,(V^{\otimes (\kn-1)}, \C) \,.
\]
On the complex $\ko C^\kn(A/\C,\C) = \gHom_\C(V^{\otimes \kn},\C)$
the differential is trivial, so 
$\gHom\,(V^{\otimes \kn}, \C)=HH^\kn(A/\C,\C)$. We finally
obtain
\[
HH^\kn(A/\C,A)= HH^\kn(A/\C,\C)\otimes V\Big/\delta_* \,HH^{\kn-1}(A/\C,\C)\;, 
\]
where the map $\delta_*$ is injective.
For the Harrison cohomology, one has to add again the condition of shuffle
invariance:
\[
\Harr^\kn(A/\C,A)= 
\Harr^\kn(A/\C,\C)\otimes V\Big/\delta_* \,\Harr^{\kn-1}(A/\C,\C)
\,.
\]
\par

{\bf Proposition:} (\cite{SchSt})
{\em Identifying the Hochschild cohomology $HH^\kb (A/\C,\C)$ with 
the tensor algebra $T V^*$ on the dual vector space  $V^*$, 
the Harrison cohomology
$\Harr^\kb (A/\C,\C)\subseteq T V^*$ consists of the primitive
elements in $T V^*$. They form a free graded Lie algebra $L$ on
$V^*$ with $V^*$ sitting in degree $-1$.
}

{\bf Proof:}
The tensor algebra $TV^*$ is a Hopf algebra with comultiplication
\[
\Delta(x_1\otimes \cdots\otimes x_\kn)
:=\sum_p
\sum_{(p,\kn-p)-{\rm shuffles}\ \sigma\hspace{-2em}}\hspace{-0.1em}
   \mbox{\rm sgn}(\sigma)\,
   (x_{\sigma(1)}\otimes\cdots\otimes x_{\sigma(p)})\otimes
   (x_{\sigma(p+1)}\otimes\cdots\otimes x_{\sigma(\kn)}) \,.
\]
It is the dual of the Hopf algebra $T^cV$ with shuffle multiplication
\[
(v_1\otimes \cdots\otimes v_p)*(v_{p+1}\otimes \cdots\otimes v_\kn)=
\sum_{(p,\kn-p)-{\rm shuffles}\ \sigma\hspace{-2em}}\hspace{-0.1em}
   \mbox{\rm sgn}(\sigma)\cdot
   v_{\sigma(1)}\otimes\cdots\otimes v_{\sigma(\kn)}\,.
\]
In particular, for any $f\in TV^*$ and $a,b \in T^cV$ one has 
$(\Delta f) (a,b) = f(a*b)$.
Hence, the condition that $f$ vanishes on shuffles is equivalent 
to $\,\Delta f = f\otimes1+1\otimes f$, i.e.\ to $f$ 
being primitive in $TV^*$.
\qed

The dimension of $\Harr^\kn(A/\C,\C)$ follows now from the dimension of
the space of homogeneous elements in the free Lie algebra, which 
was first computed in the graded case in \cite{Ree}. 

{\bf Lemma:}
{\em
$\dim_\C \Harr^\kn(A/\C,\C) = 
\frac 1\kn \sum_{d|\kn}(-1)^{\kn+\frac \kn d}\mu(d)\,m^{\frac \kn d}\,$
with $\mu$ denoting the M\"obius function.
}

{\bf Proof:}
In the free Lie algebra $L$ on $V^*$,
we choose an ordered basis $p_i$ of the even
degree homogeneous parts $L_{2\kb}$ as well as an 
ordered basis $q_i$ of the odd degree  ones.
Since $TV^*$ is the universal enveloping algebra of  $L$,
a basis for $TV^*$ 
is  given by the elements of the form $p_1^{r_1}p_2^{r_2}\cdots
p_k^{r_k}q_1^{s_1}\cdots q_l^{s_l}$ with $r_i\geq0$ and $s_i=0,1$.
In particular, if $c_\kn:=\dim L_{-\kn}$, then the Poincar\'{e}
series of the tensor algebra
\[
\sum_\kn \dim T^\kn V^*\cdot t^\kn 
=\sum_\kn m^\kn \,t^\kn={1\over 1-mt}
\]
may be alternatively described as
\[
\prod_{\kn\ {\rm even}}(1+t^\kn+t^{2\kn}+\cdots)^{c_\kn}
        \prod_{\kn\ {\rm odd}}(1+t^\kn)^{c_\kn}\,.
\]
Replacing $t$ by $-t$ and taking logarithms, the comparison of both
expressions yields
\[
-\log (1+mt)\;=\;-\sum_{\kn\ {\rm even}}c_\kn\log(1-t^\kn)+
   \sum_{\kn\ {\rm odd}}c_\kn\log(1-t^\kn)
\;=\; -\sum_{\kn}(-1)^\kn c_\kn\log(1-t^\kn)\,.
\]
Hence
\[
\sum_\kn {1\over \kn}\,(-m)^\kn\,t^\kn = 
\sum_{d, \nu}(-1)^d {1\over \nu}\,c_d\,t^{d\,\nu}\,,
\]
and by comparing the coefficients we find
\[
(-m)^\kn=\sum_{d|\kn}(-1)^d\,d\,c_d \,.
\]
Now the result follows via M\"obius inversion.
\qed
\par

We collect the dimensions
in the Poincar\'e series
\[
Q_{Z_m}(t):=\sum_{\kn\geq 1} \dim\, \Harr^{\kn}(\C\oplus V/\C,\C)\cdot t^\kn
=\sum_{\kn\geq 1}
c_n t^n\;.
\]

%%%%%%%%%%
%
%  Harrison cohomology via Noether normalization
%
%%%%%%%%%
\section{Harrison cohomology via Noether normalisation}\label{noeth}

%%%%%%%%%%%%%%%
% (noeth.norm)
%%%%%%%%%%%%%%%

\neu{noeth-norm}
Let $Y$ be a Cohen-Macaulay singularity of dimension $N$ and multiplicity
$d$; denote by $A$ its local ring. 
Choosing a {\em Noether normalisation}, i.e.\
a flat map $Y\to \C^N$ of degree $d$,
provides a regular local ring $P$ of dimension $N$ 
and a homomorphism $P\to A$ turning $A$ into
a free $P$-module of rank $d$.
Strictly speaking, this might only be possible after 
passing to an \'etale covering. Alternatively one can work in the
analytic category, see \cite{Pal} for the definition of analytic Harrison
cohomology.
\par

{\bf Proposition:}
{\em
Let $A$ be a free $P$-module as above.
If \/$M$ is any $A$-module, then
$T^\kn(A/\C,M) \cong T^\kn(A/P,M)$ for $\kn\geq 2$. 
Moreover, the latter equals
$\Harr^{\kn+1}(A/P,M)$.
}
\par

{\bf Proof:}
The Zariski-Jacobi sequence from \zitat{cotan}{gen} 
for $\C\to P \to A$ reads
\[
\cdots \longrightarrow T^\kn(A/P,M) \longrightarrow T^\kn(A/\C,M)
\longrightarrow T^\kn(P/\C,M)\longrightarrow T^{\kn+1}(A/P,M)
\longrightarrow \cdots
\]
As $P$ is regular, we have $T^\kn(P/\C,M)=0$ for $\kn\geq1$
for all $P$-modules.
On the other hand, since $A$ is flat over $P$, we may use 
\zitat{cotan}{compare}.
\qed
\par

%%%%%%%%
% (noeth.V)
%%%%%%%

\neu{noeth-V}
A rational surface singularity has minimal 
multiplicity, in the sense  
that $\mbox{\rm embdim}\,Y=\mbox{\rm mult}\,Y+\dim Y -1$.
In this situation  we may choose coordinates 
$(z_1,\dots,z_{d+1})$  such that the 
projection on the $(z_d,z_{d+1})$-plane is a Noether normalisation.
Using the above language, this means that 
$P=\C[z_d,z_{d+1}]_{(z_d,z_{d+1})}$, and
$\{1,z_1,\dots,z_{d-1}\}$ provides a basis of $A$ as a $P$-module.
\par

More generally, for a Cohen-Macaulay singularity of minimal
multiplicity we may take coordinates $(z_1,\dots,z_{d+N-1})$
such that projection on the last $N$ coordinates $(z_d,\dots,z_{d+N-1})$
is a Noether normalisation.

{\bf Lemma:}
{\em
$\;\m_A^2\;\subseteq\; \m_P\cdot\m_A$ 
\hspace{0.5em} and \hspace{0.5em}
$(\m_P\cdot\m_A)\cap P\;\subseteq \;\m_P^2\,$.
}
\par

{\bf Proof:} 
Every product $z_iz_j\in\m_A^2$ may be decomposed as
$z_iz_j=p_0 + \sum_{v=1}^{d-1}p_vz_v$ with some $p_v\in P$. Since
$\{z_1,\dots,z_{d+N-1}\}$ is a basis of 
$\raisebox{.5ex}{$\m_A$}\big/\raisebox{-.5ex}{$\m_A^2$}$, we obtain
$p_0\in\m_P^2$ and $p_v\in\m_P$ for $v\geq 1$.
The second inclusion follows from the fact that $z_1,\dots,z_{d-1}\in A$ are
linearly independent over $P$.
\qed

{\bf Proposition:}
{\em For a Cohen-Macaulay singularity of minimal multiplicity $d$
one has for $\kn\geq 1$
$T^\kn(A/\C,\C)= T^\kn(A/P,\C)=T^\kn(\C\oplus V/\C,\C)$ 
with
$V:=\raisebox{.5ex}{$\m_A$}\big/\raisebox{-.5ex}{$\m_PA$}$ 
being the $(d-1)$-dimensional vector space spanned by $z_1,\dots,z_{d-1}$.
}

{\bf Proof:} 
The equality $T^\kn(A/\C,\C)= T^\kn(A/P,\C)$ was already the subject
of Proposition \zitat{noeth}{norm} with $M:=\C$; 
it remains to treat the missing case of $\kn=1$. 
Using again the Zariski-Jacobi sequence, we have to show that
$T^0(A/\C,\C)\to T^0(P/\C,\C)$ is surjective. However, since this map is dual
to the homomorphism
$\raisebox{.5ex}{$\m_P$}\big/\raisebox{-.5ex}{$\m_P^2$}\to
\raisebox{.5ex}{$\m_A$}\big/\raisebox{-.5ex}{$\m_A^2$}$,
which is injective by the lemma above, we are done.
\vspace{0.5ex}\\
The second equality 
$T^\kn(A/P,\C)=T^\kn(\C\oplus V/\C,\C)$ follows by base change,
cf.\ \zitat{cotan}{gen}:
\[
\cdmatrix{
A&\longrightarrow &\hspace{-0.3em}\C\oplus V \cr
\mapup{{\rm flat}}&&\mapup{}\cr
P&\longrightarrow &\hspace{-0.3em}\C \cr}
\vspace{-3ex}
\]
\qed
\par

%%%%%%%%
% (noeth.rnc)
%%%%%%%

\neu{noeth-rnc}
The previous proposition reduces the cotangent cohomology with 
$\C$-coefficients
of rational surface singularities of multiplicity $d$ to that of the
fat point $Z_m$ with $m=d-1$, discussed in \zitat{cotan}{example}.
\par

{\bf Example:}
Denote by $Y_d$ the cone over the rational normal curve of degree $d$.
It may be described by the equations encoded in the condition
\[
\rank\;
\pmatrix{z_0 & z_1 &  \dots & z_{d-2}  & z_{d-1} \cr
         z_1 & z_2 &  \dots & z_{d-1}  & z_d     \cr}
\;\leq 1\;.
\]
As Noether normalisation we take the projection on the $(z_0,z_d)$-plane.
With $\deg z_i:=[i,1]\in\Z^2$, the local ring $A_d$ of $Y_d$
admits a $\Z^2$-grading. We would like to show how this grading affects the
modules
$T^\kb(A_d/\C,\C)=T^\kb(\C\oplus V/\C,\C)$ (excluding $T^0$),
i.e.\ we are going to determine the dimensions
$\, \dim T^\kb(\C\oplus V/\C,\C)(-R)$ for $R\in\Z^2$.
\hspace{1ex}\\
We know that for every $\kn$
\[
T^{\kn-1}(\C\oplus V/\C,\C)(-R)
=\Harr^{\kn}(\C\oplus V/\C,\C)(-R)\subseteq T^{\kn} V^\ast(-R)=0
\] 
unless $\kn=\kht(R)$, where $\kht(R):=R_2$ denotes 
the part carrying the standard $\Z$-grading.
Hence, we just need to calculate the numbers
\[
c_{R}:= \dim\, \Harr^{\sht(R)}(\C\oplus V/\C,\C)(-R)
\]
and can proceed as in the proof of Proposition \zitat{cotan}{example}.
Via the formal power series 
\[
\sum_{R\in\Z^2} \dim\, T^{\sht(R)}V^\ast(-R) \cdot x^R \in \C[|\Z^2|]
\]
we obtain the equation
\[
-\log \big(1+x^{[1,1]}+\cdots +x^{[d-1,1]}\big)=
   -\sum_{R\in\Z^2}(-1)^{\sht(R)}\,c_R\cdot\log(1-x^R)\,.
\]
In particular, if $\kht(R)=\kn$, then the coefficient of $x^R$ in
\[
(-1)^\kn \big(x^{[1,1]}+\cdots +x^{[d-1,1]}\big)^\kn=
(-1)^\kn \left(
\frac{x^{[d,1]}-x^{[1,1]}}{x^{[1,0]}-1}
\right)^\kn
\]
equals
$\sum_{R^\prime|R} (-1)^{\sht(R^\prime)}\,\kht(R^\prime)\cdot c_{R^\prime}$.
\vspace{0.5ex}
Again, we have to use M\"obius inversion to obtain an explicit formula for
the dimensions $c_R$.
\par

{\bf Remarks:}
\kitem
\item[(1)]
The multigraded Poincar\'{e} series
\[
\widetilde{Q}_{Z_{d-1}}(x):=\sum_{R\in\Z^2}
\dim\, \Harr^{\sht(R)}(\C\oplus V/\C,\C)(-R)\cdot x^R
= \sum_{R} c_R\,x^R
\]
is contained in the completion of the semigroup ring
$\C\big[\Z_{\geq 0}\,[1,1]+\Z_{\geq 0}\,[d-1,1]\big]$.
\item[(2)]
The cohomology groups 
$\,\Harr^{\kn}(A_d/\C,\C)(-R)$ vanish unless $\kn=\kht(R)$,
even for $n=1$. The corresponding Poincar\'{e} series 
$\widetilde Q_{Y_d}(x)$ equals
$\widetilde{Q}_{Z_{d-1}}(x) +x^{[0,1]} + x^{[d,1]}$.  
The two additional terms arise
from $\,\Harr^{1}(P/\C,\C)=T^0(P/\C,\C)$ in the exact sequence of 
\zitat{noeth}{norm}.
\kenditem
\par

%%%%%%%%
% (noeth.zero)
%%%%%%%

\neu{noeth-zero}
Let $Y$ be a Cohen-Macaulay singularity of minimal multiplicity $d\geq 3$.
\par

{\bf Lemma:}
{\em
The natural map $T^\kn(A/P,A)\to T^\kn(A/P,\C)$ is the zero map.
}
\par

{\bf Proof:}
We compute $T^\kn(A/P,\kb)$ with the reduced Harrison complex which
sits in the reduced Hochschild complex.
Using the notation of the beginning of \zitat{noeth}{V}, a reduced
Hochschild $(\kn+1)$-cocycle $f$ is, by $P$-linearity, determined
by its values on the $(\kn+1)$-tuples of the coordinates $z_1$, \dots,
$z_{d-1}$.
Suppose $f(z_{i_0}, \dots, z_{i_\kn}) \notin \m_A$. 
Since $d\geq 3$, we may choose a $z_k$ with
$k\in\{1,\dots,d-1\}$ and $k\neq i_0$. Hence
\[
\displaylines{\qquad
0=(\delta f)(z_{i_0}, \dots, z_{i_\kn},z_k)=
z_{i_0}f(z_{i_1}, \dots, z_k) \pm f(z_{i_0}, \dots, z_{i_\kn})z_k
\hfill\cr\hfill{}+
 \mbox{ terms containing products $z_iz_j$ as arguments}\;.
\qquad}
\]
Since $\m_A^2\;\subseteq\; \m_P\cdot\m_A$ by Lemma \zitat{noeth}{V}, we may
again apply $P$-linearity to see that the latter terms are contained in
$\m_P\cdot A$. Hence, modulo $\,\m_P=\m_P+\m_A^2$, these terms vanish, but
the resulting equation inside 
$V=\raisebox{.5ex}{$\m_A$}\big/\raisebox{-.5ex}{$\m_PA$}$ contradicts the
fact that $z_{i_0}$ and $z_k$ are linearly independent.
\qed
\par

{\bf Corollary:} 
{\em
\kitem
\item[{\rm(1)}]
The map $T^\kn(A/P,\m_A)\to T^\kn(A/P, A)$ is surjective. In particular,
every element of the group $T^\kn(A/P, A)$ may be represented by a cocycle
$f\colon A^{\otimes(\kn+1)}\to\m_A$.
\item[{\rm(2)}]
If $P\to A$ is $\Z$-graded with $\deg\,z_i=1$ for every $i$ 
\/{\rm(}such as for the
cone over the rational normal curve presented in Example 
{\rm\zitat{noeth}{rnc})},
then $T^\kn(A/P, A)$ sits in degree $\geq -\kn$.
\kenditem
}

%%%%%%%%%%
%
%  The cone over the rational normal curve
%
%%%%%%%%%
\section{The cone over the rational normal curve}\label{cone}

%%%%%%%%
% (cone.poinc)
%%%%%%%

\neu{cone-poinc} 
Let $Y_d$ be the cone over the rational normal curve of degree $d\geq 3$.
In Example \zitat{noeth}{rnc} we have calculated 
the multigraded Poincar\'{e} series
$\widetilde{Q}_{Y_d}(x) = 
\sum_{R} \dim\, \Harr^{\sht(R)}(A_d/\C,\C)(-R)\cdot x^R$.
The usual Poincar\'{e} series 
$Q_{Y_d}(t)=\sum_{\kn\geq 1} \dim\, \Harr^{\kn}(A_d/\C,\C)\cdot t^\kn$
is related to it
via the substitution $x^R\mapsto t^{\sht(R)}$.\\
The goal of the present section is to obtain information about
\[
P_{Y_d}(t):=\sum_{\kn\geq 1} \dim\, T^{\kn}(A_d/\C,A_d)(-R)\cdot t^\kn
\vspace{-1ex}
\]
or its multi graded version
\[
\widetilde{P}_{Y_d}(x,t):=\sum_{\kn\geq 1}\sum_{R\in\Z^2}
\dim\, T^{\kn}(A_d/\C,A_d)(-R)\cdot x^R\, t^\kn
\in \C[|\Z^2|][|t|]\,.
\]
The first series may be obtained from the latter by substituting 
$1$ for all  monomials
$x^R$, i.e.~$P_{Y_d}(t)=\widetilde{P}_{Y_d}(1,t)$.

%%%%%%%%
% (cone.toric)
%%%%%%%

\neu{cone-toric}
In \cite{AQ}, Proposition (5.2),  combinatorial
formulas have been obtained for the dimension of the vector spaces
$T^\kn(-R):=T^{\kn}(A_d/\C,A_d)(-R)$.
The point is that $Y_d$ equals the affine toric variety $Y_\sigma$ with
$\sigma$ the plane polyhedral cone
\[
\sigma:= \R_{\geq 0}\cdot (1,0) + \R_{\geq 0}\cdot (-1,d)
= \big\{ (x,y)\in\R^2\,\big|\; y\geq 0\,;\; d\,x+y\geq 0\big\}
\subseteq \R^2\,.
\]
The lattice containing the multidegrees $R$ may be identified with the
dual of the lattice $\Z^2$ inside $\R^2$,
and the results of \cite{AQ} for this special cone may be described as
follows:
\kitem
\item[(0)]
$T^0(-R)$ is two-dimensional if $R\leq 0$ on $\sigma$. It has dimension $1$
if $R$ is still non-positive on one of the $\sigma$-generators $(1,0)$ or
$(-1,d)$, but yields exactly $1$ at the other one. $T^0(-R)$ vanishes in every
other case.
\item[(1)]
$T^1(-R)$ is one-dimensional for $R=[1,1]$ and $R=[d-1,1]$;
it is two-dimensional for the degrees in between, i.e.\ for
$R=[2,1],\dots,[d-2,1]$. Altogether this means that
$\,\dim T^1=2d-4$.
\item[(2)]
The vector space $T^2$ lives exclusively in the degrees of height two.
More detailed, we have
\[
\renewcommand{\arraystretch}{1.1}
\dim T^2(-R)=\left\{\begin{array}{cl}
k-2& \mbox{for } \,R=[k,2] \;\mbox{with }\,2\leq k\leq d-1\\
d-3 &\mbox{for } \,R=[d,2] \\
2d-k-2 &\mbox{for } \,R=[k,2] \;\mbox{with } \,d+1\leq k\leq 2d-2\,.
\end{array}\right.
\vspace{-2ex}
\]
\kenditem
To formulate the result for the higher cohomology groups, we need some
additional notation.
If $R\in\Z^2$, then let $K_R$ be the finite set
\[
K_R:= \big\{ r\in\Z^2\setminus\{0\}\,\big|\;
r\geq 0 \mbox{ on } \sigma, \mbox{ but }
\,r<R \mbox{ on } \sigma\setminus\{0\} \big\}\,.
\]
Every such set $K\subseteq\Z^2$ gives rise to a complex
$C^\kb(K)$ with
\[
C^{\kn}(K):=
\Big\{\varphi\colon
\big\{(\lambda_1,\dots,\lambda_\kn)\in K^\kn\,
\big|\;\mbox{$\sum_v$}\lambda_v\in K\big\}
\to \C\,\Big|\;
\varphi \mbox{ is shuffle invariant}\Big\}\;,
\]
equipped with the modified, inhomogeneous Hochschild differential 
$d\colon C^{\kn}(K) \to C^{\kn+1}(K)$ given by
\[
\displaylines{\qquad
(d\varphi)(\lambda_0,\dots,\lambda_\kn):={}\hfill\cr\hfill
\varphi(\lambda_1,\dots,\lambda_\kn) +
\sum_{v=1}^\kn (-1)^v 
\varphi(\lambda_0,\dots,\lambda_{v-1}+\lambda_v,\dots,\lambda_\kn)
+ (-1)^{\kn+1}\varphi(\lambda_0,\dots,\lambda_{\kn-1})\,.
\qquad}
\]
Denoting the cohomology of $C^\kb(K)$ by $\HA^\kb(K)$, 
we may complete our list with the last point
\kitem
\item[(3)]
$T^\kn(-R) = \HA^{\kn-1}(K_R)\,$ for $\,\kn\geq 3$.
\kenditem 
\par

{\bf Remark:}
The explicit description of $T^2(-R)$ does almost fit into the general context
of $n\geq 3$. The correct formula is
$T^2(-R)=\raisebox{0.3ex}{$\HA^1(K_R)$}\big/ 
\raisebox{-0.3ex}{$(\span_\C K_R)^\ast$}$.
\par

%%%%%%%%
% (cone.vanish)
%%%%%%%

\neu{cone-vanish}
The previous results on $T^\kn(-R)$ 
have two important consequences.
Let
$\Lambda:=\{R\in\Z^2\,|\; R\geq 0 \mbox{ on } \sigma\}$,
$\Lambda_+:=\Lambda\setminus\{0\}$ and 
$\innt\Lambda:=\{R\in\Z^2\,|\; R> 0 \mbox{ on } \sigma\setminus \{0\}\}$.
\par

{\bf Proposition:}
{\em
Let $\kn\geq 1$.
\kitem
\item[{\rm(1)}]
$T^\kn(-R)=0$ unless $R$ is strictly positive on $\sigma\setminus\{0\}$,
i.e.\ unless $R\in\Lambda_+$.
\item[{\rm(2)}]
$T^\kn(-R)=0$ unless $\kht(R)=\kn$. In particular, $T^\kn$ is killed
by the maximal ideal of $A_d$.
\vspace{-1ex}
\kenditem
}
\par

{\bf Proof:}
(1) If $R$ is not positive on $\sigma\setminus\{0\}$, then
$K_R=\emptyset$.
\vspace{1ex}\\
(2) If $\kn-1\geq \kht(R)$, then $C^{\kn-1}(K_R)=0$ for trivial reasons.
Hence  $T^\kn$ sits in degree $\leq (-\kn)$. But this is exactly the
opposite inequality from
Corollary \zitat{noeth}{zero}(2).
\qed
\par

In particular, we may shorten our Poincar\'{e} series to
\[
\widetilde{P}_{Y_d}(x):=\sum_{\sht(R)\geq 1}
\dim\, T^{\sht(R)}(A_d/\C,A_d)(-R)\cdot x^R
\in \hat{A}_d = \C[|\Lambda|]\subseteq \C[|\Z^2|]\,.
\]
We obtain $P_{Y_d}(t)$ from 
$\widetilde{P}_{Y_d}(x)$ via the substitution $x^R\mapsto t^{\sht(R)}$.
\par

\pagebreak[2]
%%%%%%%%
% (cone.geq3)
%%%%%%%

\neu{cone-geq3}
{\bf Lemma:}
{\em 
\kitem
\item[{\rm(1)}]
Let $R\in \Z^2$ with $\kht(R)\geq 3$. Then
\[
\displaylines{\quad
\dim \,T^{\sht(R)}(-R)= 
{\kd\sum\limits_
{\smash{r\in \mbox{\rm\kf int} \Lambda,\, R-r\in\Lambda_+\hspace{-2em}}}}
(-1)^{\sht(r)-1} 
\dim \,\Harr^{\sht(R)-\sht(r)}(A_d/\C,\C)(r-R) \hfill\cr\hfill{}
+\, (-1)^{\sht(R)-1} \dim \,\HA^1(K_R) \,.\qquad}
\]
\item[{\rm(2)}]
For $R\in\Z^2$ with $\kht(R)=1$ or $2$, the right hand side 
of the above formula always yields zero.
\kenditem
}

{\bf Proof:}
(1) The vanishing of $T_Y^\kn(-R)$ for $\kn\neq\kht(R)$ together with
the equality $T_Y^\kn(-R)=\HA^{\kn-1}(K_R)$ for $\kn\geq 3$ implies that the
complex $C^\kb(K_R)$ is exact up to the first and the
$(\kht(R)-1)$-th place. In particular, we obtain
\[
\dim \,T^{\sht(R)}(-R)=\sum_{\kn\geq 1} (-1)^{\sht(R)-1+\kn} 
\dim \,C^\kn(K_R)
+ (-1)^{\sht(R)-1} \dim\, \HA^1(K_R) 
\]
where the sum is a finite one because  $C^{\geq\sht(R)}(K_R)=0$.
Now the trick is to replace the differential of the inhomogeneous complex
$C^\kb(K_R)$ by its homogeneous part $d^\prime\colon 
C^\kn(K_R)\to C^{\kn+1}(K_R)$
defined as
\[
(d^\prime\varphi)(\lambda_0,\dots,\lambda_\kn):=
\sum_{v=1}^\kn (-1)^v
\varphi(\lambda_0,\dots,\lambda_{v-1}+\lambda_v,\dots,\lambda_\kn)\,.
\]
Then $\big(C^\kb(K_R),\,d^\prime\big)$ splits into a direct sum
$\oplus_{r\in K_R} V^\kb(-r)$ with
\[
V^{\kn}(-r):=
\Big\{\varphi\colon
\big\{(\lambda_1,\dots,\lambda_\kn)\in \Lambda_+^\kn\,
\big|\;\mbox{$\sum_v$}\lambda_v=r\big\}
\to \C\,\Big|\;
\varphi \mbox{ is shuffle invariant}\Big\}\,.
\]
On the other hand, since $A_d=\C[\Lambda]$, we recognise this exactly
as the reduced complex computing $\Harr^\kb (A_d/\C,\C)(-r)$.
Hence,
\[
\dim \,T^{\sht(R)}(-R) = \sum_{\kn\geq 1,\,r\in K_R\hspace{-1em}}
(-1)^{\sht(R)-1+\kn} \dim\,\Harr^\kn(A_d/\C,\C)(-r)
\,+\, (-1)^{\sht(R)-1} \dim\, \HA^1(K_R).
\]
Finally, we replace $r$ by $R-r$ and recall  
% from Example \zitat{noeth}{rnc} 
that $\Harr^\kn(A_d/\C,\C)(-r)=0$ unless $n=\kht(r)$.
\par

(2) If $\kht(R)=2$, then the right hand side equals
$\#K_R-\#K_R=0$. If $\kht(R)=1$, then no summand at all survives.
\qed
\par

%%%%%%%%
% (cone.result)
%%%%%%%

\neu{cone-result}
Let $\,F(x):=\sum_{v=1}^{d-1} x^{[v,1]}-x^{[d,2]}=
\frac{x^{[d,1]}-x^{[1,1]}}{x^{[1,0]}-1} -x^{[d,2]}$. 
Using the Poincar\'{e} series $\widetilde{Q}_{Y_d}(x)$ of \zitat{noeth}{rnc},
we obtain the following formula:
\par

{\bf Theorem:}
{\em
The multigraded Poincar\'{e} series of the cone over the rational normal
curve of degree $d$ equals
\[
\widetilde{P}_{Y_d}(x)
\;=\;
\frac{F(x)\cdot(\widetilde{Q}_{Y_d}(x)+2)}
{\big(x^{[0,1]}+1\big)\big(x^{[d,1]}+1\big)}
-
\frac{x^{[1,1]}}{x^{[0,1]}+1}
-
\frac{x^{[d-1,1]}}{x^{[d,1]}+1}\,.
\]
}
\par

{\bf Proof:}
The previous lemma implies that
\[
\displaylines{\qquad
\widetilde{P}_{Y_d}(x) =
\sum_{\sht(R)=1,2\hspace{-1em}}\dim\,T^{\sht(R)}(-R)\cdot x^R
+ \sum_{ R\in\Lambda_+}(-1)^{\sht(R)-1} \dim \,\HA^1(K_R)\cdot x^R
\hfill\cr\hfill{}+ 
\sum_{\hspace{-2.0em}\stackrel
{\ks R\in\Lambda_+}
{\ks r\in\mbox{\rm\kf int}\Lambda,\,R-r\in\Lambda_+}
\hspace{-2.0em}}
(-1)^{\sht(r)-1}
\dim \,\Harr^{\sht(R-r)}(A_d/\C,\C)(r-R) \cdot x^R\,.
\qquad}
\]
Using the description of $T^1(-R)$ and  $T^2(-R)$
from \zitat{cone}{toric}, including
the remark at the very end,
we obtain for the first two summands
\[
\Big(2\sum_{v=1}^{d-1}x^{[v,1]} - x^{[1,1]}-x^{[d-1,1]}\Big)
+
\sum_{\ks R\in\Lambda_+}(-1)^{\sht(R)-1} \dim \,\span(K_R)\cdot x^R\,,
\]
which is equal to
\[
2\hspace{-0.5em}\sum_{R\in{\rm int}\Lambda} (-1)^{\sht(R)-1} x^R
\;+\;
\sum_{k\geq 1} (-1)^k x^{[1,k]}
\;+\;
\sum_{k\geq 1} (-1)^k x^{[kd-1,k]}\,.
\]
The third summand in the above formula for $\widetilde{P}_{Y_d}(x)$ may be
approached by summing over $r$ first. Then, substituting
$s:=R-r\in\Lambda_+$ and splitting $x^R$ into the product 
$x^r\cdot x^s$, we see that this summand is nothing else than
\[
\Big(\hspace{-0.5em}\sum_{r\in{\rm int}\Lambda} (-1)^{\sht(r)-1} x^r
\Big)\cdot
\widetilde{Q}_{Y_d}(x)\,.
\]
In particular, we obtain
\[
\widetilde{P}_{Y_d}(x)=
\Big(\hspace{-0.5em}\sum_{R\in{\rm int}\Lambda} (-1)^{\sht(R)-1} x^R
\Big)\cdot
\Big(\widetilde{Q}_{Y_d}(x)+2\Big)
\;+\;
\sum_{k\geq 1} (-1)^k x^{[1,k]}
\;+\;
\sum_{k\geq 1} (-1)^k x^{[kd-1,k]}\,.
\]
Finally, we should calculate the infinite sums. The latter two are geometric
series; they yield
$-x^{[1,1]}\big/\big(x^{[0,1]}+1\big)$ and
$-x^{[d-1,1]}\big/\big(x^{[d,1]}+1\big)$, respectively.
With the first sum we proceed as follows:
\[
\renewcommand{\arraystretch}{1.8}
\begin{array}{rcl}
\sum_{R}(-1)^{\sht(R)-1}x^R
&=& -\sum_{k\geq 1}\sum_{v=1}^{kd-1} (-1)^k x^{[v,k]} \\
&=& -\sum_{k\geq 1} (-1)^k x^{[1,k]} \,\big(x^{[kd-1,0]}-1\big)\big/
    \big(x^{[1,0]}-1\big)\\
&=& -\Big(\sum_{k\geq 1} (-1)^k x^{[kd,k]} - 
     \sum_{k\geq 1} (-1)^k x^{[1,k]}\Big)
    \Big/\Big(x^{[1,0]}-1\Big) \\
&=& \Big(x^{[d,1]}\big/\big(1+x^{[d,1]}\big) 
    - x^{[1,1]}\big/\big(1+x^{[0,1]}\big)\Big)
    \Big/\Big(x^{[1,0]}-1\Big) \\
&=& \Big( x^{[d,1]} - x^{[d+1,2]} + x^{[d,2]} - x^{[1,1]}\Big)
    \Big/\Big( \big(x^{[d,1]}+1\big)\big(x^{[1,0]}-1\big)
    \big(x^{[0,1]}+1\big)\Big)\\
% &=& \Big( \sum_{v=1}^{d-1} x^{[v,1]} -x^{[d,2]}\Big)
%     \Big/\Big( \big(x^{[d,1]}+1\big)\big(x^{[0,1]}+1\big)\Big)\\
&=& F(x)
    \Big/\Big( \big(x^{[d,1]}+1\big)\big(x^{[0,1]}+1\big)\Big)\,.
\end{array}
\vspace{-2ex}
\]
\qed
\par

%%%%%%%%
% (cone.example)
%%%%%%%

\neu{cone-example}
As example we determine the $\kht=3$ part of $\widetilde{P}_{Y_d}(x)$.
We need the first terms of $\widetilde{Q}_{Z_{d-1}}(x)$. 
By  \zitat{noeth}{rnc}
the $\kht=1$ part is just $x^{[1,1]}+\cdots +x^{[d-1,1]}$, whereas
the $\kht=2$ part is 
\[
\frac12 \big((x^{[1,1]}+\cdots +x^{[d-1,1]})^2
-(x^{[2,2]}+x^{[4,2]}+\cdots +x^{[2(d-1),2]})\big)=
\frac{(x^{[d,1]}-x^{[1,1]})(x^{[d,1]}-x^{[2,1]})}
{(x^{[1,0]}-1)^2(x^{[1,0]}+1)}
\;.
\]
Inserting this in the formula for $\widetilde{P}_{Y_d}(x)$
we finally find the grading of $T^3_{Y_d}$:
\[
\frac{(x^{[d,1]}-x^{[1,1]})(x^{[d-1,1]}-x^{[2,1]})
(x^{[d-2,1]}-x^{[2,1]})}
{(x^{[1,0]}-1)^2(x^{[2,0]}-1)}
\;.
\]
For even $d$ we get the symmetric formula
\[
(x^{[1,1]}+\cdots +x^{[d-1,1]})
(x^{[2,1]}+x^{[3,1]}+\cdots +x^{[d-3,1]}+x^{[d-2,1]})
(x^{[2,1]}+x^{[4,1]}+\cdots +x^{[d-4,1]}+x^{[d-2,1]})
\;.
\]

%%%%%%%%
% (cone.forget)
%%%%%%%

\neu{cone-forget}
Applying the ring homomorphism $x^R\mapsto t^{\sht(R)}$ to
the formula of Theorem \zitat{cone}{result} yields:
  
{\bf Corollary:}
{\em
The ordinary Poincar\'{e} series $P_{Y_d}(t)$ of the cone over 
the rational normal curve equals
\[
P_{Y_d}(t)\,=\;\Big(Q_{Y_d}(t)+2\Big)\cdot
\frac{(d-1)\,t-t^2}{(t+1)^2}
\;-\; \frac{2\,t}{t+1}\,.
\]
}
\par

%%%%%%%%%%
%
%  Hyperplane sections
%
%%%%%%%%%
\section{Hyperplane sections}\label{hps}

Choosing a Noether normalisation of an $N$-dimensional singularity
$Y$ means writing  $Y$ as 
total space of an $N$-parameter family. In this situation we have compared 
the cohomology of $Y$ with that of the $0$-dimensional special fibre.
Cutting down the dimension step by step leads to the comparison 
of the cohomology of a singularity and its hyperplane section.
\par

%%%%%%%%
% (hps.setup)
%%%%%%%

\neu{hps-setup}
First we recall the main points from \cite{BC}. Let 
$f\colon Y\to\C$ be a flat map
such that both $Y$ and the special fibre $H$ have isolated singularities.
By $T^\kn_Y$ and $T^\kn_H$ we simply denote the cotangent cohomology
$T^\kn(\CO_Y/\C,\CO_Y)$ and $T^\kn(\CO_H/\C,\CO_H)$, respectively.
\par

{\bf Main Lemma:} (\cite{BC}, (1.3))
{\em
There is a long exact sequence
\[
T^1_H \longrightarrow
T^2_Y \stackrel{\cdot f}{\longrightarrow}
T^2_Y \longrightarrow
T^2_H \longrightarrow
T^3_Y \stackrel{\cdot f}{\longrightarrow}
T^3_Y \longrightarrow
\dots\,.
\]
Moreover, $\dim\,\raisebox{0.4ex}{$T^2_Y$}\big/
            \raisebox{-0.4ex}{$f\cdot T^2_Y$}= \tau_H-e_{H,Y}$
with $\tau_H:=\dim\,T^1_H$ and $e_{H,Y}$ 
denoting the dimension of the smoothing
component containing $f$ inside the versal base space of $H$.
}
\par

This lemma will be an important tool for the comparison of the 
Poincar\'{e} series $P_Y(t)$ and $P_H(t)$ of $Y$ and $H$, respectively. 
However, since we are not only interested in the dimension, but also
in the number of generators of the cohomology groups, we introduce the
following notation.
If $M$ is a module over a local ring $(A,\m_A)$, then 
\[
\cg(M):=\dim_\C \raisebox{0.5ex}{$M$}\big/
                \raisebox{-0.5ex}{$\m_AM$}
\]
is the number of elements in a minimal generator set of $M$. 
By $P^\scg_Y(t)$
and $P^\scg_H(t)$ we denote the Poincar\'{e} series using ``$\cg$''
instead of ``$\dim$''. Similarly, $\tau_\kb^\scg:=\cg(T^1_\kb)$.
\par

{\bf Proposition:} 
{\em
\kitem
\item[\rm{(1)}]
Assume that $f\cdot T^\kn_Y=0$ for $\kn\geq 2$.
Then 
$\,P_H(t)=(1+1/t)\,P_Y(t) - \tau_Y\,(t+1) + e_{H,Y}\,t$.
\item[{\rm(2)}]
If\/ $\m_H\cdot T^\kn_H=0$ for $\kn\geq 2$, then
$\,P^\scg_H(t)=(1+1/t)\,P^\scg_Y(t) - \tau^\scg_Y\,(t+1) + 
(\tau_H^\scg-\tau_H+e_{H,Y})\,t$.
\kenditem
}
\par

{\bf Proof:}
In the first case the long exact sequence of the Main Lemma splits into
short exact sequences
\[
0\longrightarrow T^\kn_Y \longrightarrow T^\kn_H
\longrightarrow T^{\kn+1}_Y \longrightarrow 0
\]
for $\kn\geq 2$. Moreover, the assumption that $f$ annihilates
$T_Y^2$ implies that $e_{H,Y}=\tau_H-\dim\,T^2_Y$.\\
For the second part we follow the arguments
of \cite{BC}, (5.1). The short sequences have to be replaced by
\[
0\longrightarrow \raisebox{0.5ex}{$T^\kn_Y$}\big/
                 \raisebox{-0.5ex}{$f\cdot T^\kn_Y$}
\longrightarrow T^\kn_H \longrightarrow 
\ker\big[f\colon T^{\kn+1}_Y \to T^{\kn+1}_Y\big]\longrightarrow 0\,.
\]
Since $T^{\kn+1}_Y$ is finite-dimensional, the dimensions of
$\ker\big[f\colon T^{\kn+1}_Y \to T^{\kn+1}_Y\big]$ and
$\raisebox{0.5ex}{$T^\kn_Y$}\big/ \raisebox{-0.5ex}{$f\cdot T^\kn_Y$}$
are equal. Now, the claim follows from the fact that
$\raisebox{0.5ex}{$T^\kn_Y$}\big/ \raisebox{-0.5ex}{$f\cdot T^\kn_Y$}
=\raisebox{0.5ex}{$T^\kn_Y$}\big/ \raisebox{-0.5ex}{$\m_Y\, T^\kn_Y$}$,
which is a direct consequence of the assumption $\m_H\cdot T^\kn_H=0$.
\qed
\par

%%%%%%%%
% (hps.partition)
%%%%%%%

\neu{hps-partition}
We would like to apply the previous formulas to 
partition curves $H(d_1,\dots,d_r)$. They 
are defined as the wedge of the monomial curves $H(d_i)$ described by the 
equations
\[
\rank\;
\pmatrix{z_1 & z_2 &  \dots & z_{d_i-1}  & z_{d_i} \cr
         z_1 & z_3 &  \dots & z_{d_i}  & z_1^2     \cr}
\;\leq 1
\]
(\cite{BC}, 3.2).
The point making partition curves so exciting is that they occur as 
the general
hypersurface sections of rational surface singularities. Moreover, they sit
right in between the cone over the rational normal curve
$Y_d$ and the fat point $Z_{d-1}$ with $d:=d_1+\dots+d_r$. 
\par

{\bf Theorem:}
{\em
Let $H:=H(d_1,\dots,d_r)$ be a partition curve. For $\kn\geq 2$
the modules $T^\kn_H$ are annihilated by the maximal ideal $\m_H$. The
corresponding Poincar\'{e} series is
\[
P_H(t)\;=\;
\frac{\kd d-1-t}{\kd t+1} \,Q_{Z_{d-1}}(t) +\tau_H\,t -(d-1)^2\,t\,.
\vspace{-2ex}
\]
}
\par

{\bf Proof:}
We write $Y:=Y_d$ and $Z:=Z_{d-1}$.
The idea is to compare $P_H(t)$ and $P^\scg_H(t)$ which can be calculated
from $P_{Y}(t)$ and $P^\scg_{Z}(t)$, respectively.
Firstly, since $\m_YT^\kn_Y=0$ for $\kn\geq 1$, we obtain from
Proposition \zitat{hps}{setup}(1) and
Corollary \zitat{cone}{forget} that
\[
\renewcommand{\arraystretch}{1.6}
\begin{array}{rcl}
P_H(t) &=& (1+1/t)\,P_Y(t) - \tau_Y\,(t+1) + e_{H,Y}\,t\\
&=& \Big(Q_Y(t)+2\Big)
\frac{\kd (d-1)-t}{\kd t+1} -2 - (2d-4)(t+1) + e_{H,Y}\,t \\
&=& \frac{\kd d-1-t}{\kd t+1} \,Q_Z(t) -2\,t\,(d-1) + e_{H,Y}\,t\,,
\end{array}
\]
where we used that $Q_Y(t)=Q_Z(t)+2\,t$.
On the other hand, since $\m_ZT^\kn_Z=0$ for all $\kn$, we can use
the second part of Proposition \zitat{hps}{setup} to get
\[
P_Z(t)\;=\; P^\scg_Z(t)\;=\;
(1+1/t)\,P^\scg_H(t) - \tau^\scg_H\,(t+1) + e_{Z,H}\,t\,.
\]
The calculations of \zitat{cotan}{example} give  us  $P_Z(t)$ explicitly:
we have $\dim T^\kn_Z = (d-1) c_{\kn+1}-c_\kn$ and 
$\dim T^0_Z = (d-1)^2$. Therefore 
\[
P_Z(t) 
\;=\; \frac{\kd d-1-t}{\kd t}\, Q_Z(t) - (d-1)^2\,.
\]
Hence,
\[
P^\scg_H(t)=\frac{\kd d-1-t}{\kd t+1} \,Q_Z(t)
+ \tau_H^\scg\,t - \frac{t}{t+1} \Big((d-1)^2+e_{Z,H}\,t\Big)\,.
\]
Finally, we use that $\,\tau_H-e_{H,Y}=(d-1)(d-3)$ and 
$\,e_{Z,H}=(d-1)^2$ (see \cite{BC}, (4.5) and (6.3.2) respectively).
This implies the $P_H(t)$-formula of the theorem as well as 
\[
P_H(t)-P^\scg_H(t)\;=\; \big(\tau_H-\tau_H^\scg\big)\,t\,.
\]
In particular, if $\kn\geq 2$, then the modules $T^\kn_H$ have  
as dimension the number of generators, i.e.\ they are killed by the
maximal ideal.
\qed
\par

{\bf Corollary:}
{\em
The number of generators of $T^{\geq 2}$ is the same for all
rational surface singularities with fixed multiplicity $d$.
}
\par

{\bf Proof:}
Apply again Proposition \zitat{hps}{setup}(2).
\qed
\par

%%%%%%%%
% (hps.app)
%%%%%%%

\neu{hps-app}
We have seen that $\dim\,T^\kn= \cg\,T^\kn$ ($n\geq2$) for the
cone over the rational normal curve. This property holds 
for a larger class of singularities, including quotient
singularities.

{\bf Theorem:}
{\em
Let $Y$ be a rational surface singularity such that the 
projectivised tangent cone has only hypersurface singularities.
Then
the dimension of $T^\kn$ for $\kn\geq 3$ equals the number of
generators.
}
\par

{\bf Proof:}
Under the assumptions of the theorem the tangent cone
$\ko Y$ of $Y$ has also finite-dimensional $T^\kn$, $\kn\geq2$.
With $d:=\mbox{mult}(Y)$
we shall show that $\dim\,T^\kn_{\ko Y} = \dim\,T^\kn_{Y_d}$ 
for $\kn\geq 3$. 
As $Y$ is a deformation of its tangent cone $\ko Y$, 
semi-continuity implies that $\dim T^\kn_Y =
\dim T^\kn_{Y_d}$, which equals the number of generators of $T^\kn_Y$.
The advantage of working with $\ko Y$ is that it is a homogeneous
singularity, so Corollary \zitat{noeth}{zero}(2) applies.

The general hyperplane section $\ko H$ of $\ko Y$ is in general 
a non-reduced curve. In fact, it is a wedge of curves 
described by the  equations
\[
\rank\;
\pmatrix{z_1 & z_2 &  \dots & z_{d_i-1}  & z_{d_i} \cr
         z_1 & z_3 &  \dots & z_{d_i}  & 0    \cr}
\;\leq 1\,,
\]
which  is  the tangent  cone to the curve  $H(d_i)$. 
The curve $\ko H$ is also a special section of the cone over
the rational normal curve; to see this it suffices to take the cone
over a suitable divisor of degree $d$ on $\P^1$. 
Applying the Main Lemma to $\ko H$ and $Y_d$ we obtain the
short exact sequences
\[
0\longrightarrow T^\kn_{Y_d} \longrightarrow T^\kn_{\ko H}
\longrightarrow T^{\kn+1}_{Y_d} \longrightarrow 0\;,
\]
which show that for $n\geq2$ 
the dimension of $T^\kn_{\ko H}$ is the same as that of a
partition curve of multiplicity $d$. Moreover, as
the module $T^\kn_{Y_d}$ is concentrated in degree $-\kn$,
$T^{\kn+1}_{Y_d}$ in degree $-(\kn+1)$ and the  connecting homomorphism,
being induced by a  coboundary map, has degree $-1$, it
follows that  $T^\kn_{\ko H}$ 
is concentrated in degree $-\kn$. 

We now look again at $\ko H$ as hyperplane section of $\ko Y$.
The short exact sequence corresponding to the second one
in the proof of Proposition \zitat{hps}{setup},
yields that 
$\ker\big[f\colon T^{\kn+1}_{\ko Y} \to T^{\kn+1}_{\ko Y}\big]$ 
is concentrated 
in degree $-(\kn+1)$ for $n\geq2$.
The part of highest degree in $T^{\kn+1}_{\ko Y}$ is contained
in this kernel, as multiplication by $f$ increases the degree.
On the other hand, $T^{\kn+1}_{\ko Y}$ sits in degree
$\geq -(n+1)$ by Corollary \zitat{noeth}{zero}(2). Therefore
$T^{\kn+1}_{\ko Y}$ is concentrated in degree $\geq -(n+1)$
and its dimension equals 
the number of generators, which is the same as for all rational
surface singularities of multiplicity $d$. 
\qed

Note that we cannot conclude anything in the case $\kn=2$
and in fact the result does not hold for the famous 
counterexample (\cite{BC} 5.5).
\par

{\bf Corollary:}
{\em
For quotient singularities the dimension of $T^\kn$, $\kn\geq 2$,
depends only on the multiplicity. In particular, the 
Poincar\'{e} series is 
\[
P(t)\,=\;\Big(Q_{Y_d}(t)+2\Big)\cdot
\frac{(d-1)\,t-t^2}{(t+1)^2}
\;-\; \frac{2\,t}{t+1}-(\tau-2d+4)t\,.
\]
}
\par

{\bf Proof:}
For $\kn=2$ this is \cite{BC}, (Theorem 5.1.1.(3)). 
If $\kn\geq 3$, then we use the
previous theorem.
In the formula of Corollary \zitat{cone}{forget} we have then
only to introduce a correction term for $\tau=\dim\, T^1$.
\qed
\par

\newpage

%%%%%%%%%%%%%%%
%
% BIBLIOGRAPHY
%
%%%%%%%%%%%%%%%%%%
{\makeatletter
\labelsep=0pt
\def\@biblabel#1{#1}
\def\@lbibitem[#1]#2{\item[]\if@filesw
      {\let\protect\noexpand
       \immediate
       \write\@auxout{\string\bibcite{#2}{#1}}}\fi\ignorespaces}
\makeatother

}

\vfill

{\small
\parbox{9cm}{
Klaus Altmann\\
Institut f\"ur reine Mathematik der\\
Humboldt-Universit\"at zu Berlin\\
Ziegelstr.~13A\\
D-10099 Berlin, Germany\\
e-mail: altmann@mathematik.hu-berlin.de}%
\setbox0\hbox{e-mail: stevens@math.chalmers.se}\hfill\parbox{\wd0}{
Jan Stevens\\
Matematik\\
G\"oteborgs universitet\\
Chalmers tekniska h\"ogskola\\
SE-412 96 G\"oteborg, Sweden\\
e-mail: stevens@math.chalmers.se}}

\end{document}